\documentclass[12pt]{article}

\newcommand{\qed}{{\hfill\rule{4pt}{7pt}}}

\usepackage{amsmath, epsfig, cite}
\usepackage{amssymb}
\usepackage{amsfonts, color}
\usepackage{latexsym}
\usepackage{amsthm}
\newtheorem{thm}{Theorem}[section]

\newtheorem{cor}[thm]{Corollary}

\newtheorem{lem}[thm]{Lemma}

\newcommand{\pf}{\noindent{\it Proof.} }

\numberwithin{equation}{section}

\begin{document}

\begin{center}
{\Large\bf Some further $q$-series identities related to divisor functions}
\end{center}
\vskip 2mm
\centerline{\large Victor J. W. Guo$^1$ and Cai Zhang$^2$}

\vskip 3mm
\centerline{\footnotesize Department of Mathematics, East China Normal University, Shanghai 200062,
 People's Republic of China}

\vskip 1mm
\centerline
{\footnotesize
{$^1$\tt jwguo@math.ecnu.edu.cn,\quad http://math.ecnu.edu.cn/\textasciitilde{jwguo}},\qquad
 $^2${\tt 51060601083@alumni.ecnu.cn} }

\vskip 0.5cm \noindent{\small{\bf Abstract.} We give new generalizations of some $q$-series identities of Dilcher and
Prodinger related to divisor functions. Some interesting special cases are also deduced, including
an identity related to overpartitions studied by Corteel and Lovejoy.   }

\vskip 2mm
\noindent{\small{\it Keywords.} Prodinger's identity; Dilcher's identity; Lagrange interpolation; l'H\^opital's rule.

\vskip 2mm
\noindent{\it MR Subject Classifications}: Primary 11B65; Secondary 30E05   }

\section{Introduction}
In the paper \cite{Uchimura81}, Uchimura proved the following identity
\begin{align}
\sum_{k=1}^{\infty}(-1)^{k-1}\frac{q^{k+1\choose 2}}{(q;q)_k(1-q^k)}
=\sum_{k=1}^\infty \frac{q^k}{1-q^k},  \label{eq:u81}
\end{align}
where $(x;q)_N=(1-x)(1-xq)\cdots(1-xq^{N-1})$ for $N\geq 0$.
It was pointed out by Dilcher \cite{Dilcher} that \eqref{eq:u81} was known much earlier; see
\cite{Kluyver}. The identity \eqref{eq:u81} has caught the interests of several authors.
Van Hamme \cite{Hamme} (see also \cite{AU,Hoffman}) gave the following finite form:
\begin{align}
\sum_{k=1}^{n}(-1)^{k-1}{n\brack k}\frac{q^{k+1\choose 2}}{1-q^k}
=\sum_{k=1}^n \frac{q^k}{1-q^k},  \label{eq:vh}
\end{align}
where the $q$-binomial coefficient is defined as
\begin{align*}
{n\brack k}
=\begin{cases}
\displaystyle\frac{(q;q)_n}{(q;q)_k(q;q)_{n-k}},&\text{if $0\leq k\leq n$,}\\[5pt]
0,&\text{otherwise.}
\end{cases}
\end{align*}
Uchimura \cite{Uchimura} obtained a generalization of \eqref{eq:vh} as follows:
\begin{align}
\sum_{k=1}^{n}(-1)^{k-1}{n\brack k}\frac{q^{k+1\choose 2}}{1-q^{k+m}}
=\sum_{k=1}^n \frac{q^k}{1-q^k}{k+m\brack m}^{-1},\quad m\geq 0.  \label{eq:uch}
\end{align}
Dilcher \cite{Dilcher} established the following multiple series generalization of \eqref{eq:vh}:
\begin{align}
\sum_{k=1}^{n}(-1)^{k-1}{n\brack k}\frac{q^{{k\choose 2}+km}}{(1-q^k)^m}
=\sum_{1\leq k_1\leq\cdots\leq k_m\leq n} \frac{q^{k_1+\cdots+k_m}}{(1-q^{k_1})\cdots(1-q^{k_m})}.  \label{eq:dilch}
\end{align}
Prodinger \cite{Prodinger} proved that
\begin{align}
\sum_{\substack{k=0\\ k\neq m}}^{n}(-1)^{k-1}{n\brack k}\frac{q^{k+1\choose 2}}{1-q^{k-m}}
=(-1)^m q^{m+1\choose 2}{n\brack m}\sum_{\substack{k=0\\ k\neq m}}^{n} \frac{q^k}{1-q^{k-m}},\quad 0\leq m\leq n.
\label{eq:prodinger}
\end{align}
Fu and Lascoux \cite{FL,FL2} gave some further generalizations of \eqref{eq:uch}--\eqref{eq:prodinger}.
Prodinger \cite{Prodinger2} and Zeng \cite{Zeng} gave different proofs of Fu and Lascoux's identities.

In this paper, we shall give some other generalizations of Prodinger's identity
\eqref{eq:prodinger} and Dilcher's identity \eqref{eq:dilch} as well as a symmetric generalization
of \eqref{eq:vh}.

\begin{thm}\label{thm:ourdivnew} For $n\geq 0$ and $0\leq l,m\leq n$, there holds
\begin{align}
&\hskip -3mm\sum_{\substack{k=0\\ k\neq m}}^{n}{n\brack k}\frac{(q/z;q)_k
(zq^{-l};q)_{n-k}}{1-q^{k-m}} z^k   \nonumber \\
&=(-1)^m q^{m+1\choose 2}{n\brack m}(zq^{-l};q)_l(zq^{-m};q)_{n-l}\left(\sum_{k=0}^{n-l-1}
\frac{zq^{k-m} }{1-zq^{k-m}} - \sum_{\substack{k=0\\ k\neq m}}^{n}
\frac{q^{k-m}}{1-q^{k-m}} \right).  \label{eq:ourdivnew}
\end{align}
\end{thm}

Noticing that
$$
\lim_{z\to 0}(q/z;q)_k z^k=(-1)^k q^{k+1\choose 2},
$$
when $l=0$ and $z$ tends to $0$, the identity \eqref{eq:ourdivnew}
reduces to \eqref{eq:prodinger}.

\begin{thm}\label{thm:dilch}
For $m,n\geq 1$, there holds
\begin{align}
&\sum_{k=1}^{n}{n\brack k}\frac{(q^m/z;q)_k(z;q)_{n-k}}{(zq^{-m};q)_{m+n}(1-q^k)^{m}}z^k  \nonumber\\
&=-\sum_{k_1=1}^n\frac{q^{k_1}}{(1-zq^{k_1-1})(1-q^{k_1})}
\sum_{k_2=1}^{k_1}\frac{q^{k_2}}{(1-zq^{k_2-2})(1-q^{k_2})}
\cdots
\sum_{k_m=1}^{k_{m-1}}\frac{q^{k_m}}{(1-zq^{k_m-m})(1-q^{k_m})}.  \label{eq:multi}
\end{align}
\end{thm}

Similarly, when $z$ tends to $0$, the identity \eqref{eq:multi}
reduces to \eqref{eq:dilch}. It should be mentioned here that \eqref{eq:multi}
is not a consequence of the following key identity appearing in \cite[(2)]{Zeng}:
$$
\sum_{1\leq i_1\leq i_2\leq\cdots\leq i_m\leq N}a_{i_1}a_{i_2}\cdots a_{i_m}
=\sum_{k=1}^N \left(\prod_{\substack{j=1,j\neq k}}^N (1-a_j/a_k)^{-1}\right)a_k^m,
$$
for the right-hand side of \eqref{eq:multi} cannot be written as a complete
symmetric function of $m$ indeterminates $a_1,a_2,\ldots,a_m$.

Our third theorem is the following symmetric generalization of the $l=m=0$
case of Theorem~\ref{thm:ourdivnew}.
\begin{thm}\label{thm:sym}
For $m,n\geq 0$, there holds
\begin{align}
&\hskip -3mm  \sum_{k=1}^{n}\frac{(q/z;q)_{k}(vq^m;q)_{k}(z;q)_{n-k}(z;q)_m}
{(q;q)_k (v;q)_{k} (q;q)_{n-k}(q^k;q)_{m+1}}z^k
-\sum_{k=1}^{m}\frac{(q/z;q)_{k} (vq^n;q)_{k}(z;q)_{m-k}(z;q)_n}
{(q;q)_k (v;q)_{k}(q;q)_{m-k}(q^k;q)_{n+1}}z^k \nonumber\\
&=\frac{(1-zq^{-1})(z;q)_m(z;q)_n}{(q;q)_{m}(q;q)_{n}}
\left(\sum_{k=1}^m\frac{q^k}{(1-zq^{k-1})(1-q^k)}-\sum_{k=1}^n\frac{q^k}{(1-zq^{k-1})(1-q^k)}\right).
\label{eq:symsym}
\end{align}
\end{thm}

The rest of the paper is organized as follows. We shall prove Theorem \ref{thm:ourdivnew}
in Section 2 and give some special cases in Section 3. Theorem \ref{thm:dilch}
will be proved in Section 4, and a class of binomial sums will be evaluated in Section 5.
Finally, the proof of Theorem \ref{thm:sym} and some interesting conclusions will be given in Section 6.

\section{Proof of Theorem \ref{thm:ourdivnew}}
Recall that the Lagrange interpolation formula states that a polynomial of degree $\leq n$
that passes through the $n+1$ points $(x_0,f(x_0)),(x_1,f(x_1)),\ldots,(x_n,f(x_n))$ is given by
\begin{align*}
f(x)=\sum_{k=0}^n f(x_k)\prod_{\substack{j=0\\ j\neq k}}^n\frac{x-x_j}{x_k-x_j}
\end{align*}
(see, for example, \cite{JJ}).

For $0\leq l\leq n$, consider the Lagrange interpolation formula for $(xz; q)_{n-l}$ at the values
$q^{-i}$ ($0\leq i\leq n$) of $x$, we immediately obtain

\begin{lem}For $n\geq 0$ and $0\leq l\leq n$, there holds
\begin{align}
\sum_{k=0}^{n}{n\brack k}\frac{(q/z;q)_k (zq^{-l};q)_{n-k}}{1-xq^k}
z^k =\frac{(q;q)_n(xz;q)_{n-l}(zq^{-l};q)_l}{(x;q)_{n+1}}. \label{eq:15}
\end{align}
\end{lem}

\medskip
\noindent{\it Proof of Theorem \ref{thm:ourdivnew}.} From \eqref{eq:15} it follows that
\begin{align}
&\hskip -3mm \sum_{\substack{k=0\\ k\neq m}}^{n}{n\brack k}\frac{(q/z;q)_k (zq^{-l};q)_{n-k}}{1-xq^k} z^k  \nonumber \\
&=\frac{(q;q)_n(xz;q)_{n-l}(zq^{-l};q)_l}{(x;q)_{n+1}}-{n\brack m}\frac{(q/z;q)_m (zq^{-l};q)_{n-m}}{1-xq^m} z^m  \nonumber\\
&=\frac{1}{(x;q)_m(xq^{m+1};q)_{n-m}} \nonumber\\
&\quad{}\times\left(\frac{(q;q)_n(xz;q)_{n-l}(zq^{-l};q)_l}{1-xq^m}
-{n\brack m}\frac{(q/z;q)_m (zq^{-l};q)_{n-m}(x;q)_m(xq^{m+1};q)_{n-m}}{1-xq^m}
z^m\right). \label{eq:2term}
\end{align}
It is easy to see that
$$
(q;q)_n(xz;q)_{n-l}(zq^{-l};q)_l
-{n\brack m}(q/z;q)_m (zq^{-l};q)_{n-m}(x;q)_m(xq^{m+1};q)_{n-m}z^m
=0
$$
if $x=q^{-m}$. By l'H\^opital's rule, we have
\begin{align*}
&\hskip -3mm \lim_{x\to q^{-m}}\left(\frac{(q;q)_n(xz;q)_{n-l}(zq^{-l};q)_l}{1-xq^m}-{n\brack m}
\frac{(q/z;q)_m (zq^{-l};q)_{n-m}(x;q)_m(xq^{m+1};q)_{n-m}}{1-xq^m}z^m\right)
\nonumber \\
&=\lim_{x\to q^{-m}}(q;q)_n(zq^{-l};q)_l\sum_{k=0}^{n-l-1} \frac{zq^{k-m} (xz;q)_{n-l}}{1-xzq^k}  \\
&\quad{} -\lim_{x\to q^{-m}}{n\brack m}(q/z;q)_m
(zq^{-l};q)_{n-m}z^m \sum_{\substack{k=0\\ k\neq m}}^{n}
\frac{q^{k-m} (x;q)_m(xq^{m+1};q)_{n-m}}{1-xq^k}
\nonumber \\
&=(q;q)_n(zq^{-m};q)_{n-l}(zq^{-l};q)_l\left(\sum_{k=0}^{n-l-1}
\frac{zq^{k-m} }{1-zq^{k-m}} - \sum_{\substack{k=0\\ k\neq m}}^{n}
\frac{q^{k-m}}{1-q^{k-m}} \right).
\end{align*}
Therefore, letting $x\to q^{-m}$ in \eqref{eq:2term} and applying
the relation
$$
\frac{(q;q)_n}{(q^{-m};q)_m (q;q)_{n-m}}=(-1)^m q^{m+1\choose
2}{n\brack m},
$$
we complete the proof.  \qed

\section{Consequences of Theorem \ref{thm:ourdivnew} }
Letting $l=0$ in Theorem \ref{thm:ourdivnew}, we get
\begin{cor}\label{thm:ourdiv}
For $n\geq 0$ and $0\leq m\leq n$, there holds
\begin{align}
\sum_{\substack{k=0\\ k\neq m}}^{n}{n\brack k}\frac{(q/z;q)_k (z;q)_{n-k}}{1-q^{k-m}} z^k
=(-1)^m q^{m+1\choose 2}{n\brack m}(zq^{-m};q)_n \left(\sum_{k=0}^{n-1} \frac{zq^{k-m} }{1-zq^{k-m}}
- \sum_{\substack{k=0\\ k\neq m}}^{n} \frac{q^{k-m}}{1-q^{k-m}} \right).  \label{eq:ourdiv}
\end{align}
\end{cor}

Letting $n\to\infty$ in \eqref{eq:ourdiv}\,, we obtain
\begin{align}
\sum_{\substack{k=0\\ k\neq m}}^{\infty}\frac{(q/z;q)_k z^k}{(q;q)_k(1-q^{k-m})}
=(-1)^m q^{m+1\choose 2}\frac{(zq^{-m};q)_m}{(q;q)_m}
\left(\sum_{k=0}^{\infty} \frac{zq^{k-m} }{1-zq^{k-m}}
- \sum_{\substack{k=0\\ k\neq m}}^{\infty} \frac{q^{k-m}}{1-q^{k-m}} \right).  \label{eq:ourinfty}
\end{align}
For $m=0$, the above identity \eqref{eq:ourinfty} reduces to
\begin{align}
\sum_{k=1}^{\infty}\frac{(q/z;q)_k z^k}{(q;q)_k(1-q^{k})}
=\sum_{k=0}^{\infty} \frac{zq^{k} }{1-zq^{k}}
- \sum_{k=1}^{\infty} \frac{q^{k}}{1-q^{k}}.  \label{eq:ourinfty2}
\end{align}
Now, letting $z=-q$ in \eqref{eq:ourinfty2}\,, we are led to
\begin{align*}
\sum_{k=1}^{\infty}\frac{(-1;q)_{k} (-q)^k}{(q;q)_k(1-q^{k})}
=-\sum_{k=1}^{\infty} \frac{2q^{k}}{1-q^{2k}},
\end{align*}
of which a combinatorial interpretation was given by Corteel and Lovejoy \cite[Theorem 4.4]{CL}.

Letting $m=0$ in Theorem \ref{thm:ourdivnew}, we get

\begin{cor} For $n\geq 0$ and $0\leq l\leq n$, there holds
\begin{align}
\sum_{k=1}^{n}{n\brack k}\frac{(q/z;q)_k (zq^{-l};q)_{n-k}}{1-q^k} z^k
=(zq^{-l};q)_n\left(\sum_{k=1}^{n-l}
\frac{zq^{k-1}}{1-zq^{k-1}}-\sum_{k=1}^n\frac{q^k}{1-q^k}\right).   \label{eq:cor-our}
\end{align}
\end{cor}

Furthermore, letting $l=0$ in \eqref{eq:cor-our}\,, and observing that
\begin{align}
\frac{zq^{k-1}}{1-zq^{k-1}}-\frac{q^k}{1-q^k}
=-\frac{q^{k}(1-zq^{-1})}{(1-zq^{k-1})(1-q^{k})}, \label{eq:lab}
\end{align}
we obtain
\begin{cor} For $n\geq 0$, there holds
\begin{align}
\sum_{k=1}^{n}{n\brack k}\frac{(q/z;q)_k (z;q)_{n-k}}{1-q^k} z^k
=(zq^{-1};q)_{n+1}\sum_{k=1}^n\frac{q^{k}}{(1-zq^{k-1})(1-q^{k})}.  \label{eq:cor1}
\end{align}
\end{cor}

\section{Proof of Theorem \ref{thm:dilch}}

\pf We proceed by induction on $m$ and on $n$. For $m=1$, the identity \eqref{eq:multi} reduces to
\eqref{eq:cor1}. Assume that \eqref{eq:multi} holds for some $m\geq 1$. We need to show that
it also holds for $m+1$, namely,
\begin{align}
&\sum_{k=1}^{n}{n\brack k}\frac{(q^{m+1}/z;q)_k(z;q)_{n-k}}{(zq^{-m-1};q)_{m+n+1}(1-q^k)^{m+1}}z^k  \nonumber\\
&=-\sum_{k_0=1}^n\frac{q^{k_0}}{(1-zq^{k_0-1})(1-q^{k_0})}
\sum_{k_1=1}^{k_0}\frac{q^{k_1}}{(1-zq^{k_1-2})(1-q^{k_1})}
\cdots
\sum_{k_m=1}^{k_{m-1}}\frac{q^{k_m}}{(1-zq^{k_m-m-1})(1-q^{k_m})}.  \label{eq:multi2}
\end{align}
We shall prove this induction step \eqref{eq:multi2} by induction on $n$, following the proofs
in \cite{Hoffman} and \cite[Theorem 4]{Dilcher}.

For $n=1$, both sides of \eqref{eq:multi2} are equal to
$$
-\frac{q^{m+1}}{(zq^{-m};q)_{m+1}(1-q)^{m+1}},
$$
and it is true. We assume that \eqref{eq:multi2} holds for $n-1$. In order to show that it also holds for $n$,
we have to check that the difference between \eqref{eq:multi2} for $n$ and \eqref{eq:multi2} for $n-1$
is a true identity. This difference is
\begin{align}
&\sum_{k=1}^{n}\frac{(q^{m+1}/z;q)_k}{(1-q^k)^{m+1}}z^k
\left({n\brack k}\frac{(z;q)_{n-k}}{(zq^{-m-1};q)_{m+n+1}}-{n-1\brack k}\frac{(z;q)_{n-k-1}}{(zq^{-m-1};q)_{m+n}}\right)  \nonumber\\
&=-\frac{q^{n}}{(1-zq^{n-1})(1-q^{n})}
\sum_{k_1=1}^{n}\frac{q^{k_1}}{(1-zq^{k_1-2})(1-q^{k_1})}
\cdots
\sum_{k_m=1}^{k_{m-1}}\frac{q^{k_m}}{(1-zq^{k_m-m-1})(1-q^{k_m})}.  \label{eq:diff}
\end{align}
Now using the relation
\begin{align*}
{n\brack k}\frac{(z;q)_{n-k}}{(zq^{-m-1};q)_{m+n+1}}-{n-1\brack k}\frac{(z;q)_{n-k-1}}{(zq^{-m-1};q)_{m+n}}
={n\brack k}\frac{(zq^{-1};q)_{n-k}(1-q^k)q^{n-k}}{(zq^{-m-1};q)_{m+n}(1-q^n)(1-zq^{n-1})},
\end{align*}
we see that \eqref{eq:diff} is equivalent to \eqref{eq:multi} with $z$ being replaced by $zq^{-1}$.
This proves that \eqref{eq:multi2} holds for all $n\geq 1$, and consequently \eqref{eq:multi} holds for all $m\geq 1$
and all $n\geq 1$.     \qed

\medskip
Similarly to the proof of Theorem \ref{thm:dilch}, we can prove that
\begin{thm}\label{thm:dilch2}
For $m,n\geq 1$, there holds
\begin{align}
&\sum_{k=1}^n{n\brack k}\frac{(q^{m}/z;q)_k(z/q;q)_{n-k}}{(zq^{-m};q)_{n+m-1}(1-q^k)^m}z^k \nonumber\\
&=\sum_{k_1=1}^{n}\frac{q^{k_1}}{(1-zq^{k_1-2})(1-q^{k_1})}
\sum_{k_2=1}^{k_1}\frac{q^{k_2}}{(1-zq^{k_2-3})(1-q^{k_2})} \cdots
\left(\sum_{k_m=1}^{k_{m-1}-1}\frac{zq^{k_m-m}}{1-zq^{k_m-m}}
-\sum_{k_m=1}^{k_{m-1}}\frac{q^{k_{m}}}{1-q^{k_m}}\right).  \label{eq:dilch2}
\end{align}
\end{thm}

By taking the limit as $n\to\infty$ in Theorems \ref{thm:dilch} and \ref{thm:dilch2},
we obtain the following two results.

\begin{cor} For $m\geq 1$, there holds
\begin{align*}
&\sum_{k=1}^{\infty}\frac{(q^m/z;q)_k}{(zq^{-m};q)_{m}(q;q)_k(1-q^k)^{m}}z^k  \nonumber\\
&=-\sum_{k_1=1}^\infty\frac{q^{k_1}}{(1-zq^{k_1-1})(1-q^{k_1})}
\sum_{k_2=1}^{k_1}\frac{q^{k_2}}{(1-zq^{k_2-2})(1-q^{k_2})}
\cdots
\sum_{k_m=1}^{k_{m-1}}\frac{q^{k_m}}{(1-zq^{k_m-m})(1-q^{k_m})}.
\end{align*}
\end{cor}

\begin{cor} For $m\geq 1$, there holds
\begin{align*}
&\sum_{k=1}^{\infty}\frac{(q^{m}/z;q)_k}{(zq^{-m};q)_{m-1}(q;q)_k(1-q^k)^{m}}z^k  \nonumber\\
&=\sum_{k_1=1}^\infty\frac{q^{k_1}}{(1-zq^{k_1-2})(1-q^{k_1})}
\sum_{k_2=1}^{k_1}\frac{q^{k_2}}{(1-zq^{k_2-3})(1-q^{k_2})}
\cdots
\left(\sum_{k_m=1}^{k_{m-1}-1}\frac{zq^{k_m-m}}{1-zq^{k_m-m}}
-\sum_{k_m=1}^{k_{m-1}}\frac{q^{k_{m}}}{1-q^{k_m}}\right).
\end{align*}
\end{cor}

\section{Some binomial identities}
As was pointed out by Andrews and Uchimura \cite{AU}\,, the partial sum
$$
H_n(q):=\sum_{k=1}^n\frac{q^k}{1-q^k}
$$
may be considered as a $q$-analogue of the harmonic number $H_n:=\sum_{k=1}^n 1/k$. Thus, the identity
\eqref{eq:vh} is a $q$-analogue of the celebrated identity
\begin{align}
\sum_{k=1}^{n}(-1)^{k-1}{n\choose k}\frac{1}{k}=\sum_{k=1}^n\frac{1}{k}. \label{eq:trigo}
\end{align}
(see, for example, \cite[(1.46)]{Gould}). Dilcher \cite{Dilcher} has obtained a multiple generalization
of \eqref{eq:trigo} by taking the limit as $q\to 1$ in \eqref{eq:dilch} as follows:
\begin{align}
\sum_{k=1}^{n}(-1)^{k-1}{n\choose k}\frac{1}{k^m}
=\sum_{1\leq k_m\leq k_{m-1}\leq\cdots\leq k_1\leq n}\frac{1}{k_1k_2\cdots k_m}.  \label{eq:dilch-noq}
\end{align}
To obtain a further generalization of \eqref{eq:dilch-noq}\,, we multiply both sides of
\eqref{eq:multi} by $(1-q^k)^{2m}$, replace $z$ by $q^x$, and let $q$ tend to $1$. Then we get
the following result.
\begin{cor} For $m,n\geq 1$, there holds
\begin{align}
\sum_{k=1}^{n}{n\choose k}\frac{(m-x)_k(x)_{n-k}}{(x-m)_{m+n}k^m}
=-\sum\frac{1}{k_1k_2\cdots k_m (x+k_1-1)(x+k_2-2)\cdots(x+k_m-m)},  \label{eq:multi-noq}
\end{align}
where $(x)_N=x(x+1)\cdots (x+N-1)$ and the sum ranges over all
integers $1\leq k_m\leq k_{m-1}\leq\cdots\leq k_1\leq n$.
\end{cor}
It is clear that when $x$ tends to $\infty$, the identity \eqref{eq:multi-noq} (multiplying both sides by $-x^m$)
reduces to Dilcher's identity \eqref{eq:dilch-noq}\,. On the other hand, taking the limit as $x\to m$ in
\eqref{eq:multi-noq}\,, we are led to
\begin{align*}
&\hskip -3mm \sum_{k=1}^{n}{n\choose k}\frac{(k-1)!(m+n-k-1)!}{k^{m}}  \\
&=\sum_{1\leq k_m\leq k_{m-1}\leq\cdots\leq k_1\leq n}
\frac{(m-1)!(m+n-1)!}{k_1k_2\cdots k_m (k_1+m-1)(k_2+m-2)\cdots k_m}.
\end{align*}

Similarly, we can derive the following result from \eqref{eq:dilch2}\,.
\begin{cor} For $m,n\geq 1$, there holds
\begin{align*}
&\hskip -3mm \sum_{k=1}^{n}{n\choose k}\frac{(m-x)_k(x-1)_{n-k}}{(x-m)_{m+n-1}k^m} \\
&=\sum\frac{\sum_{j=1}^{k_{m-1}-1}\frac{1}{x+j-m}-\sum_{j=1}^{k_{m-1}}\frac{1}{j}}
{k_1k_2\cdots k_{m-1} (x+k_1-2)(x+k_2-3)\cdots(x+k_{m-1}-m)},
\end{align*}
where the sum is over all integers $1\leq k_{m-1}\leq\cdots\leq k_1\leq n$.
\end{cor}

\section{Proof of Theorem \ref{thm:sym}}

We start with the following identity appearing in Guo and Zeng \cite{GZ}:
\begin{align}
\frac{(xz,yz;q)_m}{(q,xyz;q)_m }
\sum_{k=0}^{n}\frac{(x,y,vq^m;q)_k(z;q)_{n-k}}
{(q,v,xyzq^m;q)_k (q;q)_{n-k}}z^k
=\frac{(xz,yz;q)_n}{(q,xyz;q)_n }
\sum_{k=0}^{m}\frac{(x,y,vq^n;q)_k(z;q)_{m-k}} {(q,v,xyzq^n;q)_k
(q;q)_{m-k}}z^k, \label{eq:zeng}
\end{align}
which can be obtained applying Sears's transformation \cite[p.~360, (III.15)]{GR}.
Letting $y=q/z$ in \eqref{eq:zeng}, we get
\begin{equation*}
\sum_{k=0}^{n}\frac{(q/z;q)_{k}(vq^m;q)_{k}(z;q)_{n-k}(xz;q)_m}
{(q;q)_k (v;q)_k (q;q)_{n-k}(xq^k;q)_{m+1}}z^k
=\sum_{k=0}^{m}\frac{(q/z;q)_{k}(vq^n;q)_{k}(z;q)_{m-k}(xz;q)_n}
{(q;q)_k (v;q)_{k}(q;q)_{m-k}(xq^k;q)_{n+1}}z^k,
\end{equation*}
which can be rewritten as
\begin{align}
&\hskip -3mm\sum_{k=1}^{n}\frac{(q/z;q)_{k}(vq^m;q)_{k}(z;q)_{n-k}(xz;q)_m}
{(q;q)_k (v;q)_{k} (q;q)_{n-k}(xq^k;q)_{m+1}}z^k
-\sum_{k=1}^{m}\frac{(q/z;q)_{k} (vq^n;q)_{k}(z;q)_{m-k}(xz;q)_n}
{(q;q)_k (v;q)_{k}(q;q)_{m-k}(xq^k;q)_{n+1}}z^k \nonumber\\
&=\frac{(z;q)_{m}(q;q)_{n}(xz;q)_n(xq;q)_{m}-(z;q)_{n}(q;q)_{m}(xz;q)_m(xq;q)_{n}}
{(q;q)_{m}(q;q)_{n}(xq;q)_{m}(xq;q)_{n}(1-x)}.  \label{eq:newsym}
\end{align}
It is obvious that
$$
(z;q)_{m}(q;q)_{n}(xz;q)_n(xq;q)_{m}-(z;q)_{n}(q;q)_{m}(xz;q)_m(xq;q)_{n}=0\quad
\text{if $x=1$}.
$$
Let $x\to 1$ in \eqref{eq:newsym}. Then applying l'H\^opital's rule and the relation \eqref{eq:lab},
we complete the proof of Theorem \ref{thm:sym}.

Letting $v\to 0$ or $v\to\infty$ in \eqref{eq:symsym}, we obtain the following two corollaries.

\begin{cor}
For $m,n\geq 0$, there holds
\begin{align*}
&\hskip -3mm  \sum_{k=1}^{n}\frac{(q/z;q)_{k}(z;q)_{n-k}(z;q)_m}
{(q;q)_k (q;q)_{n-k}(q^k;q)_{m+1}}z^k
-\sum_{k=1}^{m}\frac{(q/z;q)_{k}(z;q)_{m-k}(z;q)_n}
{(q;q)_k (q;q)_{m-k}(q^k;q)_{n+1}}z^k \\
&=\frac{(1-zq^{-1})(z;q)_m(z;q)_n}{(q;q)_{m}(q;q)_{n}}
\left(\sum_{k=1}^m\frac{q^k}{(1-zq^{k-1})(1-q^k)}-\sum_{k=1}^n\frac{q^k}{(1-zq^{k-1})(1-q^k)}\right).
\end{align*}
\end{cor}

\begin{cor}
For $m,n\geq 0$, there holds
\begin{align*}
&\hskip -3mm  \sum_{k=1}^{n}\frac{(q/z;q)_{k}(z;q)_{n-k}(z;q)_m}
{(q;q)_k (q;q)_{n-k}(q^k;q)_{m+1}}q^{mk}z^k
-\sum_{k=1}^{m}\frac{(q/z;q)_{k}(z;q)_{m-k}(z;q)_n}
{(q;q)_k (q;q)_{m-k}(q^k;q)_{n+1}}q^{nk}z^k \\
&=\frac{(1-zq^{-1})(z;q)_m(z;q)_n}{(q;q)_{m}(q;q)_{n}}
\left(\sum_{k=1}^m\frac{q^k}{(1-zq^{k-1})(1-q^k)}-\sum_{k=1}^n\frac{q^k}{(1-zq^{k-1})(1-q^k)}\right).
\end{align*}
\end{cor}

Taking the limit as $z\to 0$, the above two corollaries reduce to
\begin{cor}For $m,n\geq 0$, there holds
\begin{align}
&\hskip -3mm  \sum_{k=1}^{n}\frac{(-1)^k q^{{k+1\choose 2}}}
{(q;q)_k (q;q)_{n-k}(q^k;q)_{m+1}}
-\sum_{k=1}^{m}\frac{(-1)^k q^{{k+1\choose 2}}}
{(q;q)_k (q;q)_{m-k}(q^k;q)_{n+1}} \nonumber\\
&=\frac{1}{(q;q)_{m}(q;q)_{n}}
\left(\sum_{k=1}^m\frac{q^k}{1-q^k}-\sum_{k=1}^n\frac{q^k}{1-q^k}\right). \label{cor:001}
\end{align}
\end{cor}

\begin{cor}For $m,n\geq 0$, there holds
\begin{align}
&\hskip -3mm  \sum_{k=1}^{n}\frac{(-1)^k q^{mk+{k+1\choose 2}}}
{(q;q)_k (q;q)_{n-k}(q^k;q)_{m+1}}
-\sum_{k=1}^{m}\frac{(-1)^k q^{nk+{k+1\choose 2}}}
{(q;q)_k (q;q)_{m-k}(q^k;q)_{n+1}} \nonumber\\
&=\frac{1}{(q;q)_{m}(q;q)_{n}}
\left(\sum_{k=1}^m\frac{q^k}{1-q^k}-\sum_{k=1}^n\frac{q^k}{1-q^k}\right).  \label{cor:002}
\end{align}
\end{cor}
It is clear that both \eqref{cor:001} and \eqref{cor:002} reduce to \eqref{eq:vh}
if $m=0$. Moreover, if we perform the substitution $q\to q^{-1}$, the identities
\eqref{cor:001} and \eqref{cor:002} can respectively be written as
\begin{align*}
&\hskip -3mm  \sum_{k=1}^{n}\frac{(-1)^k q^{{k+1\choose 2}+(m-n)k}}
{(q;q)_k (q;q)_{n-k}(q^k;q)_{m+1}}
-\sum_{k=1}^{m}\frac{(-1)^k q^{{k+1\choose 2}+(n-m)k}}
{(q;q)_k (q;q)_{m-k}(q^k;q)_{n+1}} \\
&=\frac{1}{(q;q)_{m}(q;q)_{n}}
\left(\sum_{k=1}^m\frac{1}{1-q^k}-\sum_{k=1}^n\frac{1}{1-q^k}\right), \\[5pt]
&\hskip -3mm  \sum_{k=1}^{n}\frac{(-1)^k q^{{k+1\choose 2}-nk}}
{(q;q)_k (q;q)_{n-k}(q^k;q)_{m+1}}
-\sum_{k=1}^{m}\frac{(-1)^k q^{{k+1\choose 2}-mk}}
{(q;q)_k (q;q)_{m-k}(q^k;q)_{n+1}} \\
&=\frac{1}{(q;q)_{m}(q;q)_{n}}
\left(\sum_{k=1}^m\frac{1}{1-q^k}-\sum_{k=1}^n\frac{1}{1-q^k}\right).
\end{align*}
Finally, letting $n\to\infty$ in \eqref{cor:001} and \eqref{cor:002}\,, we obtain two different
generalizations of \eqref{eq:u81}\,:
\begin{align*}
\sum_{k=1}^{\infty}\frac{(-1)^{k-1} q^{{k+1\choose 2}}}
{(q;q)_k (q^k;q)_{m+1}}
-\sum_{k=1}^{m}\frac{(-1)^{k-1} q^{{k+1\choose 2}}} {(q;q)_{m-k}(1-q^k)}
&=\frac{1}{(q;q)_{m}}\sum_{k=m+1}^\infty\frac{q^k}{1-q^k}, \\ 
\sum_{k=1}^{\infty}\frac{(-1)^{k-1} q^{mk+{k+1\choose 2}}}
{(q;q)_k (q^k;q)_{m+1}}
&=\frac{1}{(q;q)_{m}}\sum_{k=m+1}^\infty\frac{q^k}{1-q^k}.
\end{align*}

\vskip 5mm
\noindent{\bf Acknowledgements.} This work was partially
supported by the Fundamental Research Funds for the Central Universities, Shanghai Rising-Star Program (\#09QA1401700),
Shanghai Leading Academic Discipline Project (\#B407), and the National Science Foundation of China (\#10801054).

\renewcommand{\baselinestretch}{1}

\end{document}